\documentclass[12pt]{amsart}
\usepackage{amssymb}

\swapnumbers

\newtheorem*{Thmn}{Theorem}
\newtheorem*{Lemman}{Lemma}
\newtheorem{Thm}[subsection]{Theorem}
\newtheorem{Cor}[subsection]{Corollary}

\newtheorem{Prop}[subsection]{Proposition}

\theoremstyle{definition}

\newtheorem{Ex}[subsection]{Example}

\theoremstyle{remark}
\newtheorem{Remark}[subsection]{Remark}

\DeclareMathOperator{\re}{Re}
\DeclareMathOperator{\im}{Im}

\newcommand{\R}{\ensuremath{\mathbb R}}
\newcommand{\N}{\ensuremath{\mathbb N}}
\newcommand{\C}{\ensuremath{\mathbb C}}
\newcommand{\abs}[1]{\left\vert#1\right\vert}
\newcommand{\chf}{\ensuremath{\mathbf{1}}}
\newcommand{\M}{\ensuremath{\mathcal{M}}}
\newcommand{\m}{\ensuremath{\mathbf{m}}}
\newcommand{\cum}{\ensuremath{\mathbf{c}}}
\newcommand{\Rtr}{\ensuremath{\mathcal{R}}}

\renewcommand{\epsilon}{\varepsilon}

\title[Linearization of the CLT operator]
{The linearization of the central limit operator in free probability theory}

\author[M. Anshelevich]{Michael Anshelevich}

\thanks{This work is supported in part by the Fannie and John Hertz
  Foundation Fellowship}

\address{Department of Mathematics, University at California\\
Berkeley, CA 94720, USA}

\email{mashel@math.berkeley.edu}


\date{October 6, 1998}

\begin{document}

\begin{abstract}
We interpret the Central Limit Theorem as a fixed point theorem for a certain operator, and
consider the problem of linearizing this operator. In classical as well as in free
probability theory \cite{VDN92}, we consider two methods giving such a linearization, and
interpret the result as a weak form of the CLT. In the classical case the analysis involves
dilation operators; in the free case more general composition operators appear.
\end{abstract}

\maketitle

\section{Introduction}
The ubiquity of the normal distribution as indicated by the Central Limit Theorem (CLT) is
a somewhat mysterious result. One of the possible explanations for it is an interpretation
of the CLT as a fixed-point theorem;  for a plethora of approaches see
\cite{Tro59,Gol76,HW84,Bar86,Swa91,Sin92}. The starting point for this analysis is the
following weak form of the CLT for independent identically distributed random variables.
Let the operator $T$ be defined on probability measures on $\R$ by $T \mu = (\mu \ast \mu)
\circ S_{\frac{1}{\sqrt{2}}}$. Here $\ast$ is the convolution, and $S_r$ (for ``scaling'')
is the dilation operator, $d(\mu \circ S_r)(x) = d\mu(r^{-1}x)$. We call this operator $T$
the \emph{central limit} operator. The theorem follows from the CLT, and is well known.

\begin{Thmn}
The fixed points of the operator $T$ are the scaled normal distributions $\chi \circ S_t$.
If $\mu$ is a probability measure with zero mean and unit variance, then the iterations
$T^n \mu$ weakly converge to $\chi$.
\end{Thmn}

In this approach, the starting point for the study of the CLT is the investigation of the
operator $T$. This operator is clearly non-linear, and as the first approximation we
consider the properties of the linearization of this operator \cite{Sin92}. These
properties are of course well-known (although we have not found adequate references; but
see \cite{Sin92} and also \cite{Sin76}). However, there is now a different version of
probability theory, with its own CLT, which has not been investigated to the same degree.
This is the free probability theory of Voiculescu (for an introduction, see e.g.
\cite{VDN92}), which,  in particular, turns out to describe the behavior of certain large
random matrices. In this theory the notion of (commutative) independence is replaced by the
notion of free independence (for the operator-theoretic definition and the motivation
behind it see e.g. \cite{VDN92}). Now the classical convolution can be defined in terms of
independence as follows: $\mu \ast \nu$ is the distribution of the sum of two independent
random variables with distributions $\mu$ and $\nu$ (and it is a theorem that the
distribution of the sum of independent random variables depends only on the distributions
of the summands). Correspondingly, in free probability theory one defines the (additive)
free convolution of measures $\mu \boxplus \nu$ as the distribution of the sum of two
freely independent random variables with distributions $\mu$ and $\nu$
\cite{Voi85,Voi86,BV93,VDN92} (and the above comment applies). Thus the \emph{free central
limit operator} is  $T(\mu) = (\mu \boxplus \mu) \circ S_{\frac{1}{\sqrt{2}}}$. One of the
main technical differences between classical and free theories of probability is that the
operator of convolution with a given measure is linear, while the operator of taking a free
convolution with a given measure is highly non-linear. However, the above operator $T$ is
non-linear even in the classical case, and thus one can expect similarities between
linearizations of the classical and free versions of this operator.

We have two somewhat different approaches at our disposal. The original one, initiated and
largely developed by Voiculescu \cite{Voi86,BV93, VDN92} (see also \cite{Maa92}), is to
define a certain operation, called the \emph{$R$-transform}, on the space of analytic
functions, which linearizes the additive free convolution. Thus this operation is an
analogue of the logarithm of the Fourier transform in the classical case. Another
approach, due to Speicher \cite{Spe90} and developed, among others, by Speicher and Nica
\cite{Spe94,Nic95}, is to use a certain analogue of the classical combinatorial
moment-cumulant formula. This approach is somewhat less general, but the parallel with the
classical situation is more explicit.

In what follows we want to indicate the parallels between the classical and the free case.
Therefore, whenever appropriate, we will use the same notation for both cases. The
situations where (important) differences between the two theories arise will also be
indicated.

\noindent {\bf Acknowledgments.} We would like to thank Prof.\ D.-V.~Voiculescu for
suggesting the problem as well as many helpful discussions. We would also like to thank
Prof.\ N.G.~Makarov for some suggestions.

\section{The Combinatorial Approach}
\label{sec:Comb}
\subsection{Notation}
Let $\mu$ be a probability measure. We denote by $T$ both the central limit operator
$T(\mu)  = (\mu \ast \mu) \circ S_{\frac{1}{\sqrt{2}}}$ and the free central limit
operator $T(\mu) = (\mu \boxplus \mu) \circ S_{\frac{1}{\sqrt{2}}}$. It will be clear from
the context which one is meant. We denote by $\mathcal{T}$ the manifestations of $T$ on
auxiliary spaces: in Section~\ref{sec:Comb}, the spaces of sequences; in
Section~\ref{sec:Anal}, the spaces of continuous functions. Precise definitions will be
given at appropriate times. Also, we denote by $\chi$ the appropriate normal
distributions: standard Gaussian in the classical context and the standard Wigner
semicircle law \cite{Voi85} in the free context.

\subsection{Background}
For a measure $\mu$, its $n$-th moment is $m^\mu_n = \int x^n d \mu(x)$. In this section we
consider only probability measures whose moments of all orders are finite. In fact,
throughout most of the section we disregard the non-uniqueness and identify the measure
with its collection of moments. Thus let $\M$ be the space of all one-sided real-valued
sequences. We will call the elements of $\M$ the moment sequences and denote them by $\m
:= \{m_i\}_{i=1}^\infty$, even though only some of them are moment sequences of measures.
On $\M$ we set up the topology of entrywise convergence; this is the weak$^\ast$-topology
on the space $\M$ as the dual of the space of the ``eventually $0$'' sequences, and it
turns $\M$ into a topological vector space. Note that if a sequence of elements of $\M$ do
in fact correspond to measures, and if its limit corresponds to a \emph{unique} measure,
then one has weak convergence of the corresponding measures \cite{Dur91}.

For every such moment sequence $\m$ there is also the corresponding \emph{free cumulant}
sequence $\cum := \{c_i\}_{i=1}^\infty$ \cite{Spe94,Nic95} determined by
\begin{equation}
\label{eq:freemc}
\text{(free)} \qquad m_k = \sum_{\substack{\pi \in \mathcal{P}_{nc}(k)\\
 \pi = \{B_1, \ldots, B_n\}}} \prod_{j=1}^{n} c_{\abs{B_j}}
\end{equation}
Here $\mathcal{P}_{nc}(k)$ is the set of \emph{noncrossing} partitions of the set
$\{1,\ldots,k\}$ \cite{Kre72,Spe90,Spe94,Nic95}, which can be described as follows: these
are partitions of the vertices of an $k$-gon such that the vertices in each class can be
connected by lines inside the $k$-gon so that the lines for different classes do not
cross. Also, $B_i$-s denote the classes of the partition $\pi$, and $\abs{B_i}$ denotes
the number of elements of $B_i$.

The classical cumulant sequence can also be described by a similar formula \cite{Shi96,Nic95},
namely
\begin{equation}
\label{eq:classmc}
\text{(classical)} \qquad m_k = \sum_{\substack{\pi \in \mathcal{P}(k)\\
 \pi = \{B_1, \ldots,  B_n\}}} \prod_{j=1}^{n} (\abs{B_j} - 1)! \; c_{\abs{B_j}}
\end{equation}
where $\mathcal{P}(k)$ is the collection of all partitions of $\{1,\ldots,k\}$.

Let us denote the transformation from the moment sequence to the cumulant sequence
determined by the formula \eqref{eq:freemc} (resp., \eqref{eq:classmc}) by $\Rtr :
\{m_i\}_{i=1}^\infty \rightarrow \{c_i\}_{i=1}^\infty$. We call $\Rtr^{-1}$ the
\emph{cumulant-moment} transform. Note that $\Rtr$ is given implicitly; there are also
explicit formulas~\cite{Nic95}. Note also that $k$-th moment depends only on the cumulants
of orders less than or equal to $k$, and vice versa. Therefore both $\Rtr$ and $\Rtr^{-1}$
are continuous bijections $\M \rightarrow \M$; however, we will think of the domain of
$\Rtr$ as moment sequences and of its range as cumulant sequences.

The point of the transformation from the moment to the cumulant sequence is that for
sequences corresponding to probability measures, the appropriate action of the operator
$T$ on the cumulant side is linear. Indeed,
\begin{equation*}
\label{eq:cumT}
c^{\mu \boxplus \nu}_k = c^\mu_k + c^\nu_k  \qquad  \text{and}  \qquad c^{\mu \circ S_r}_k
= r^k c^\mu_k
\end{equation*}
where $c^\eta$ are the cumulants of the measure $\eta$. That is,
\begin{equation}
\label{cumulants}
c_k^{T \mu} = 2^{1 - \frac{k}{2}} c_k^\mu
\end{equation}
Thus define, on the space of cumulant sequences, the operator $\mathcal{T^R}$ by
$(\mathcal{T^R}(\cum))_k  = 2^{1 - k/2} c_k$, and on $\M$ the operator $\mathcal{T} =
\Rtr^{-1} \mathcal{\circ T^R \circ R}$. Clearly, since $\mathcal{T^R}$ is linear, in order
to linearize the operator $\mathcal{T}$, we are interested in the linearization of the
cumulant-moment transform $\Rtr^{-1}$.

In the sequel, by a \emph{linearization} of a map $A$ at a point $x$ we mean its
G\^{a}teaux derivative: $(D_x A) (y) = \lim_{\epsilon \rightarrow 0} \frac{A(x +
\epsilon y) - A(x)}{\epsilon}$ when the limit exists in the appropriate topology. Note also
that from~\eqref{cumulants}, a fixed point of $T$ for which the cumulant sequence is
defined must have all the cumulants other than the second one equal to 0. In the classical
case, this describes the Gaussian distributions; in the free case, this describes the free
normal distributions, which are the dilations of the Wigner semicircle law \cite{Spe90}.

\begin{Prop}
The linearization of the cumulant-moment transform at $($the cumulant series corresponding
to$)$ the normal distribution $\chi$ $($respectively, standard Gaussian in the classical
case and standard Wigner semicircular distribution in the free case \cite{VDN92}$)$ is
given by a $($formal$)$ infinite lower-triangular matrix $A = (a_{ij})_{i,j = 1}^\infty$,
where

\begin{enumerate}
\item
In the classical case, $a_{n + 2k, n}$ is $(n-1)!$ times the number of partitions of
$(n+2k)$ elements into classes exactly one of which contains $n$ elements and the
remaining $k$ classes are pairs.
\item
In the free case, $a_{n + 2k, n}$ is the number of noncrossing partitions of $(n + 2k)$
elements into classes exactly one of which contains $n$ elements and the remaining $k$
classes are pairs.
\end{enumerate}
In both cases $a_{ij}=0$ if $j > i$ or $(i-j)$ is odd. For explicit values, see
Theorem~\ref{thm:Herm}.
\end{Prop}

\begin{proof}
The value of the $n$-th cumulant is a polynomial function of the first $n$ moments only,
and vice versa. Thus in the topology of entrywise convergence, the differentials of both
$\Rtr$ and $\Rtr^{-1}$ exist.

Given two moment sequences $\m^o, \m^d$, we define the sequence $\{f(\m^o,
\m^d)_n\}_{n=1}^\infty$ recursively by
\begin{equation}
\label{CombDeriv}
m_k^d = \sum_{\substack{\pi \in \mathcal{P}_{nc} (k)\\
 \pi = \{B_1, \ldots, B_n\}}} \sum_{i =1}^n
\left( \prod_{j \neq i} c_{\abs{B_j}}^o \right) f(\m^o, \m^d)_{\abs{B_i}}
\end{equation}
where $\{c^o_i\} = \Rtr(\m^o)$ are the free cumulants. Then
\begin{equation*}
\sum_{\substack{\pi \in \mathcal{P}_{nc} (k)\\
 \pi = \{B_1, \ldots, B_n\}}} \prod_{i=1}^n
 \left( c_{\abs{B_i}}^o + \epsilon f(\m^o, \m^d)_{\abs{B_i}} \right)
 = m_k^o + \epsilon m_k^d + o(\epsilon)
\end{equation*}
Note that if $\m^o = \m^\mu, \m^d = \m^\nu$, then the last expression above is just
$\m^{\mu + \epsilon \nu} + o(\epsilon)$. So the sequence $\{f(\m^o,
\m^d)_n \}_{n=1}^\infty$ is the derivative of the moment-cumulant transform $\Rtr$ at $\m^o$
in the direction $\m^d$. For $\m^o = \m^\chi$, the free standard normal (semicircular)
distribution, $c^\chi_i = \delta_{i2}$, and so~\eqref{CombDeriv} becomes
\begin{equation}
m_k^d = \sum_{n=1}^k a_{k,n} f(\m^d)_n
\end{equation}
where $f(\m^d) := f(\m^\chi, \m^d)$ and $a_{n + 2k,n}$ is the number of noncrossing
partitions of $(n + 2k)$ elements into classes exactly one of which has $n$ elements and
the remaining $k$ classes are pairs. Thus $m = Af$, where $A$ is the lower-triangular
matrix $(A)_{i,j} = a_{i,j}$.

In the classical case, we start with
\begin{equation}
m_k^d = \sum_{\substack{\pi \in \mathcal{P} (k)\\
 \pi = \{B_1, \ldots, B_n\}}} (\abs{B_i}-1)! \sum_{i=1}^n
\left( \prod_{j \neq i} c_{\abs{B_j}}^o \right) f(\m^o, \m^d)_{\abs{B_i}}
\end{equation}
and by the same sort of reasoning see that the derivative of $\Rtr$ at the standard
Gaussian is given by the lower-triangular matrix $A$ with $a_{n + 2k, n} = (n-1)! \times $
(the number of partitions of $(n + 2k)$ elements into classes exactly one of which contains
$n$ elements and the remaining $k$ classes are pairs).
\end{proof}
As stated above, the operator $\mathcal{T^R}$ is linear. It is easy to see that its
spectrum is discrete. Its eigenvectors are the cumulant sequences $\xi_j =
\{\delta_{ij}\}_{i=1}^\infty$, for $j=1,2,\ldots$, with corresponding eigenvalues
$2^{1- j/2}$. Therefore for the linearization of operator $\mathcal{T}$, the eigenvalues
are the same, and the eigenvectors are the moment sequences $e_j
= \{a_{i,j}\}_{i=1}^\infty$, where $a_{i,j}$-s are defined in the above theorem. In fact,
these are true moment sequences, and so give the eigenfunctions for the central limit
operator $T$.

\begin{Thm}
\label{thm:Herm}
On the space of measures with all moments finite, the linearization of the operator $T$ has
eigenvalues $2^{1 - n/2}, n = 1, 2, \ldots$. The corresponding eigenfunctions are
absolutely continuous with respect to the Lebesgue measure, with densities:

\begin{enumerate}
\item
In the classical case, $\frac{d^n}{dx^n} e^{-x^2/2} = e^{-x^2/2} H_n(x)$, multiples of the
Hermite polynomials  \cite{Sin92}.
\item
In the free case, $\chf_{[-2,2]}(t) \frac{1}{\sqrt{4 - t^2}} T_n(t/2)$, multiples of the
Chebyshev polynomials of the first kind.
\end{enumerate}
\end{Thm}

\begin{proof}
In the classical case, $a_{n+2k, n} = (n-1)! \; \times$ the number of partitions of
$(n+2k)$ objects into one class of $n$ elements and $k$ classes of 2 elements. It is easy
to see that $a_{n+2k, n} = \frac{(n+2k)!}{n k! 2^k}$ and $a_{k,n}=0$ for $k <n$ or $(k-n)$
odd. Therefore for fixed $n$ the Fourier transform (defined by $\sum_{j=0}^\infty
\frac{1}{j!} m_j (it)^j = \int e^{itx} d \mu(x)$) of the $n$-th eigenfunction of $T$ is
\begin{equation*}
\sum_{k=0}^\infty a_{n+2k, n}
\frac{1}{(n+2k)!} (it)^{n+2k} = \frac{1}{n} (it)^n \exp(- t^2/2)
\end{equation*}
and the sum converges absolutely. Thus the eigenfunctions are the multiples of Hermite
polynomials $\frac{d^n}{dx^n} e^{-x^2/2} \,dx = e^{-x^2/2} H_n(x) \,dx$ (note that these
are not exactly what one usually means by the Hermite functions).

In the free case, $a_{n+2k, n} = $ the number of noncrossing partitions of \mbox{$(n+2k)$}
objects into one class of $n$ elements and $k$ classes of 2 elements. It has been
calculated by Kreweras \cite{Kre72} to be $a_{n+2k, n} = \binom{n+2k}{k}$ (one uses an
inductive argument based on the following fact: a partition $\pi$ with a class of $n$
elements is noncrossing iff each of the $n$ intervals in the complement of this class is a
union of complete classes of $\pi$, and $\pi$ restricted to each of these intervals is
noncrossing). The Cauchy transform (defined by $\sum_{j=0}^\infty m_j z^{-(j+1)} = \int
\frac{d \mu(x)}{z - x}$) (\cite{Akh65,VDN92}, see also the next section)  of the $n$-th
eigenfunction is $\sum_k \binom{n+2k}{k} z^{-(n+2k+1)}$. For $z \in \C^+$, the series
converges absolutely for $\abs{z} >2$. Its integral is
\begin{equation*}
F_n (z) = - \sum_k \frac{1}{n+k} \binom{n + 2k-1}{k} z^{-(n+2k)} = -\sum_k
\frac{1}{n+2k} \binom{n + 2k}{k} z^{-(n+2k)}
\end{equation*}
In particular, for $n=1$ we have $F_1(z) = - \sum_k \frac{1}{k+1} \binom{2k}{k}
z^{-(2k+1)}$. Thus ${F_1(z)}^2 = - F_1(z)z - 1$. Therefore $F_1$ is related to the
generating function for the Catalan numbers \cite{Rio68}, and is in fact $ \frac{-z +
\sqrt{z^2 - 4}}{2}$. Similarly

\begin{Lemman}
For $n \geq 1$ the integral of the Cauchy transform of the $n$-th eigenfunction is
$-\frac{1}{n} \left( \frac{z - \sqrt{z^2 - 4}}{2} \right)^n$, i.e. $ -\sum_k
\frac{1}{n+2k} \binom{n + 2k}{k} z^{-(n+2k)}= - \frac{1}{n}
\left( \frac{z - \sqrt{z^2 - 4}}{2} \right)^n$
\end{Lemman}

\begin{proof}[Proof of the Lemma]
The series converges absolutely for $\abs{z} > 2$. The proof is by induction, using the
identity $\left( \frac{z - \sqrt{z^2 - 4}}{2} \right)^n \left( \frac{z - \sqrt{z^2 -
4}}{2} \right)^m = \left( \frac{z - \sqrt{z^2 - 4}}{2} \right)^{(n + m)}$. By equating
coefficients, we have to prove the combinatorial identity
\begin{equation*}
\sum_{\substack{k+l=t \\
k,l \geq 0}} \frac{n}{n+2k} \binom{n + 2k}{k} \frac{m}{m+2l}
\binom{m + 2l}{l} = \frac{n + m}{n + m + 2t} \binom{n + m + 2t}{t}
\end{equation*}
for $m, n = 1, 2, \ldots$ and $t=0, 1, \ldots$. But this is a particular case of the
generalized Vandermonde (also known as Rothe) identity (\cite[Sec. 4.5]{Rio68}, see also
\cite{GK66}).
\renewcommand{\qed}{}
\end{proof}
The above moments determine a unique distribution \cite{Dur91}, and we can see directly
that the $n$-th eigenfunction is related to the Chebyshev polynomials of the first kind,
namely it is a scalar multiple of $\chf_{[-2,2]}(t) \frac{1}{\sqrt{4 - t^2}} T_n(t/2)
\,dt$, where $T_n(t) = \cos^{-1}(n \cos t)$.
\end{proof}

\begin{Remark}
Besides being eigenfunctions of the operator $DT$ on a topological vector space, the above
functions in fact form orthogonal bases in (smaller) Hilbert spaces. The Hermite functions
$\frac{d^n}{dx^n} e^{-x^2/2} = e^{-x^2/2} H_n(x)$ for $n=0, 1, \ldots$ form an orthogonal
basis in $L^2(e^{x^2/2} \,dx)$, while the Chebyshev functions of the first kind
$\chf_{[-2,2]}(t) \frac{1}{\sqrt{4 - t^2}} T_n(t/2)$ for $n = 0, 1, \ldots$ form an
orthogonal basis in the space $L^2(\chf_{[-2,2]}(t) \sqrt{4-t^2} \,dt)$. Thus, in the
classical case the operator $DT$ is a compact self-adjoint operator on $L^2(e^{x^2/2}
\,dx)$, while in the free case $DT$ is a compact self-adjoint operator on
$L^2(\chf_{[-2,2]}(t) \sqrt{4-t^2} \,dt)$.
\end{Remark}
We have the following interpretation of the CLT.

\begin{Cor}
We say that a measure $\mu$ is in $L^2 (\varphi)$ if $\mu \ll \varphi\,dx$ and $\frac{d
\mu}{dx} \in L^2 (\varphi)$. Note that for all the measures in $L^2(e^{x^2/2} \,dx)$ or
$L^2(\chf_{[-2,2]}(t) \sqrt{4-t^2} \,dt)$, the moments of all orders are finite.

\begin{enumerate}
\item
On the space of probability distributions in $L^2(e^{x^2/2} \,dx)$, with mean $0$ and
variance $1$, the normal distribution $\chi$ as a fixed point of the central limit operator
$T$ is strictly spectrally stable, that is, the differential of the operator at this point
has the spectrum inside the unit disc.
\item
On the space of probability distributions in $L^2(\chf_{[-2,2]}(t) \sqrt{4-t^2} \,dt)$,
with mean $0$ and variance $1$, the free normal (semicircular) distribution $\chi$ as a
fixed point of the free central limit operator $T$ is strictly spectrally stable.
\end{enumerate}
\end{Cor}

\begin{proof}
Classical case: \cite{Sin92} The Hermite functions $\frac{d^n}{dx^n} e^{-x^2/2} =
e^{-x^2/2} H_n(x)$ for $n=0, 1, \ldots$ form an orthogonal basis in $L^2(e^{x^2/2} \,dx)$.
The conditions on the moments of the distribution in the hypothesis mean that in the linear
approximation we consider only the perturbations $f$ with 0th, 1st, 2nd moments equal to 0.
This means precisely that $f$ is orthogonal to $e^{-x^2/2}$ and (the densities of) the
first two eigenfunctions of $T$ . Since the eigenvalues of $T$ are $2^{1 - n/2}$, all
eigenvalues for $n > 2$ are less than $1$.

Free case: the Chebyshev functions $\chf_{[-2,2]}(t) \frac{1}{\sqrt{4 - t^2}} T_n(t/2)$
for $n = 0, 1, \ldots$ form an orthogonal basis in the space $L^2(\chf_{[-2,2]}(t)
\sqrt{4-t^2} \,dt)$. Again the eigenvalues are less than $1$ for all but the first $3$ of
these. But here, the hypothesis correspond to the orthogonality to the first $3$ Chebyshev
functions only if all the distributions considered are \emph{supported in the same
interval} $[-2, 2]$. Note, however, that due to the results in \cite{BV95}, this
restriction is weaker than it appears.
\end{proof}

Thus we would expect that on this subspace the fixed point is attracting, just as the CLT
states.

\begin{Remark}
The  moments $a_{n+2k,k}$ can be calculated in a way similar to the above for the setting
of the $R_q$ transforms \cite{Nic95}, related to $q$-independence. However, there does not
seem to be in that case a nice recurrence formula, and so the corresponding eigenfunctions
are not calculated directly. Moreover, note that the orthogonality conditions in the above
Corollary are different in the classical and the free cases: in the free case the inner
product is given by the free normal (semicircular) distribution, while in the classical
case it is the inverse of the normal distribution. Also, we are not aware of any standard
interpolation between the Hermite polynomials and the Chebyshev polynomials of the first
kind. Thus one would not necessarily expect to have a similar construction for the
interpolations between free and classical cases, e.g. related to $q$-independence.
\end{Remark}

\section{The Analytic Approach}
\label{sec:Anal}
In this section we consider the problem of linearizing the operator $T$ by analytic means.
First we briefly go over the classical situation.

\subsection{Classical Picture}
Let $\alpha \in (0, 2]$, $\beta = 2^{- 1/\alpha}$, $T_\alpha \mu = (\mu \ast \mu) \circ
S_\beta$. For $\varphi_\alpha =$ $\alpha$-strictly stable distribution \cite{Shi96,Dur91}
(the skewness coefficient does not appear explicitly in the sequel and so is not included
in the notation), $T_\alpha (\varphi_\alpha) = \varphi_\alpha$. Then the differential of
$T_\alpha$ at $\varphi_\alpha$ is
\begin{equation*}
DT_\alpha  \nu = 2 (\nu \ast \varphi_\alpha) \circ S_\beta
\end{equation*}
Taking the Fourier transforms,
\begin{equation}
\label{Ftrans}
\widehat{DT_\alpha  \nu} (t)= 2 \hat{\nu}(\beta t) \hat{\varphi}_\alpha (\beta t)
\end{equation}
Also by stability $\hat{\varphi}_\alpha^2(\beta t) = \hat{\varphi}_\alpha(t)$. Therefore
for $\hat{\nu}(t) = h(t) \hat{\varphi}_\alpha(t)$, the right-hand-side expression
in~\eqref{Ftrans} is $2 h(\beta t) \hat{\varphi}_\alpha (t)$. For $h$, on the space of
continuous functions the eigenfunctions are $h(t) = t^a$, $a \in \C$, $\re a > 0$ or
$a=0$, with eigenvalues $2 \cdot \beta^a= 2 \cdot 2^{-a/\alpha}$, corresponding to
$\hat{\nu_a}(t) = t^a \hat{\varphi}_\alpha (t)$. Here we use the principal branch of the
logarithm. Now let $\beta=1/\sqrt{2}$, i.e. $\varphi_\alpha = \varphi_2 = \chi$. In this
case, among all the eigenfunctions we can distinguish the integer values of $a$ as follows:
among the functions $t^a$, the smooth ones are precisely those for $a \in \N$. Thus among
all $\nu_a$, the ones whose densities decay faster than any polynomial are just the
$\nu_n$-s, $n \in \N$ \cite{Shi96,Dur91}. These measures are manifestly in $L^2(e^{x^2/2}
\,dx)$. They are $\nu_n = \frac{d^n}{dx^n} e^{-x^2/2} \,dx$, with eigenvalues $2^{1 - n/2}$, and we obtain the result of the
previous section.

\subsection{Free Picture}
In the free probability picture, the main device is the $R$-transform, introduced by
Voiculescu \cite{Voi86,BV93,VDN92}. Given a measure $\mu$, for $z \in \C \backslash
\text{supp}(\mu)$ one defines the Cauchy transform (sometimes called Stieltjes or Borel
transform) of $\mu$ by $G_\mu(z) = \int \frac{d \mu(t)}{z-t}$. For positive $\mu$ the
Cauchy transform is an analytic map $\C^+ \leftrightarrow \C^-$; it has the property
$\bar{G}_\mu(z) = G_\mu (\bar{z})$. The measure can be reconstructed from its Cauchy
transform by taking the boundary values $- \frac{1}{\pi} \im G_\mu(x+0i)$~\cite[Ch.3,
Addenda and Problems]{Akh65},~\cite[3.1]{Hor90}.

On a nontangential neighborhood of 0 (Stolz angle) in $\C^-$, we can define $K_\mu (w) =
G_\mu^{-1}(w)$, and the $R$-transform $R_\mu(w) = K_\mu(w) - \frac{1}{w}$. $R_\mu$ is an
analytic map $\C^- \rightarrow \C^-$ on a nontangential neighborhood of 0. The main
property of the $R$-transform is that it also linearizes the additive free convolution:
$R_{\mu \boxplus \nu}(w) = R_\mu(w) + R_\nu(w)$. In fact, if all the moments of a measure
$\mu$ are finite, then $R_\mu(w) = \sum_{i=1}^\infty c_i^\mu w^{i-1}$, where $c_i$ are the
free cumulants. Also $R_{\mu \circ S_r}(w) = r R_\mu(r w)$. Thus the action of the operator
$T_\alpha$ on the $R$-transform side is just $R_{T_\alpha \mu}(w) = 2\beta R_\mu(\beta w)$,
and in particular it is linear. Here we again define the operator
\begin{equation*}
T_\alpha (\mu) = (\mu \boxplus \mu) \circ S_\beta
\end{equation*}
for $\alpha \in (0,2]$, and $\varphi_\alpha$ = free $\alpha$-strictly stable distribution
\cite{BV93,Pat95,BPB96}. By the above observations, the linearization of $T_\alpha$ is
again given by the linearization of the $R$-transform.

The quasilinear differential equation governing the behavior of free convolution
semigroups has first appeared in \cite[Theorem 4.3]{Voi86}. Here we need a variant of that
theorem. The proof is quite similar to \cite{Voi86}.

\begin{Thm}
Let $\nu$ be a freely infinitely divisible probability measure \cite{Voi86,VDN92}. Let
$\psi$ be a function analytic in $\C^+$. For any point $z \in \C^+$, for small enough $t$
the function $G(z,t)$, which is the functional inverse of $K_\nu (z) + t \psi(z)$, is
defined at $z$. Consequently for all $z \in \C^+$,
\begin{equation}
G_\nu'(z)\psi(G_\nu(z)) + \frac{\partial G}{\partial t} (z,0) = 0
\end{equation}
\end{Thm}

\begin{proof}
By \cite[Proposition 5.12]{BV93} for $\nu$ freely infinitely divisible, $G_\nu$ maps
$\C^+$ conformally onto a domain. Let $\Omega$ be a bounded domain in $\C^+$ whose closure
is contained in $\C^+$. Choose $t$ so that for $z \in \Omega$, $t \abs{\psi'(G_\nu(z))} <
1$. Then the function $z + t \psi(G_\nu (z))$ is univalent on $\Omega$ and invertible on
its image. Denote this image by $\Omega'$ and this inverse function by $f_t$. Since, for
$z \in \Omega$,
\begin{equation*}
(K_\nu + t \psi)(G_\nu (z)) = z + t \psi(G_\nu (z))
\end{equation*}
we also have
\begin{equation*}
(K_\nu + t \psi)(G_\nu (f_t(z)) = z
\end{equation*}
for $z \in \Omega'$. Consequently we can define $G(z,t)
= G_\nu (f_t(z))$ for $z \in
\Omega'$.

As one possible construction, define
\begin{equation*}
\Omega_n = \left\{ z \big| -n \leq \re z \leq n, \frac{1}{n} \leq \im z \leq n \right\}
\end{equation*}
Let $t_n$ satisfy the condition above and also $t_n \abs{\psi(G_\nu (z))} < 1/n$ for $z
\in \Omega_n$. Then $G(z,t)$ is defined on
\begin{equation*}
\Omega_n' = \left\{ z \big| -n+ \frac{1}{n} \leq \re z \leq n- \frac{1}{n}, \frac{2}{n}
\leq \im z \leq n- \frac{1}{n} \right\}
\end{equation*}
These domains exhaust $\C^+$. On their common domains, since their inverse is analytic in
$t$, $G(z,t)$ is differentiable in $t$; the derivative exists as a limit in the topology
of uniform convergence on compact sets in the upper half plane.

By definition $G(K_\nu + t \psi(z), t) = z$ for $z \in G_\nu(\Omega)$. Differentiating
with respect to $t$ at $t=0$, we get
\begin{equation*}
\frac{\partial G}{\partial z}(K_\nu(z), 0) \psi(z) +
\frac{\partial G}{\partial t}(K_\nu(z), 0) = 0
\end{equation*}
for $z \in G_\nu(\C^+)$. Substituting $G_\nu(z)$ for $z$, we get the required equation,
for $z \in \C^+$.
\end{proof}

\begin{Remark}
In the above theorem no consideration is given to the positivity or even existence of the
distribution corresponding to the Cauchy transform $G(z,t)$. While we do not know of a
satisfactory description of these, there are various conditions.

\subsubsection{Necessary Conditions}
\label{nec}
For $K_\nu + \psi$ to correspond to a positive measure, it is necessary that (1) the
nontangential limit of $z \psi(z)$ as $z \rightarrow 0$ be 0, and (2) there exist a Stolz
angle at 0, $\Gamma \subset \C^+$ s.t. $(K_\nu + \psi)(\Gamma) \subset \C^-$ \cite{BV93}.

\subsubsection{Sufficient Conditions}
\label{suf}

\begin{enumerate}
\item
\label{deform}
The following are the known cases where $K_\nu + \psi$ corresponds to a positive measure:
(1) $\psi$ is an $R$-transform of a positive (hence necessarily freely infinitely
divisible) measure \cite{Voi86}. (2) $\nu$ is the free normal (semicircular) distribution;
$\psi$ is analytic in a neighborhood of the unit disc, sufficiently small, and
$\psi(\bar{z}) = \bar{\psi}(z)$~\cite{BV95}.
\item
For a measure $\mu$ and a distribution $\nu$, the $R$-transform of $((1 - \epsilon)\mu +
\epsilon \nu)$ is
\begin{equation*}
R_{(1 - \epsilon)\mu + \epsilon \nu} = R_\mu(w) - \epsilon K_\mu'(w) \cdot (G_\nu -
G_\mu)(K_\mu(w)) + o(\epsilon)
\end{equation*}
Denote $\psi = K_\mu'(w) \cdot (G_\nu - G_\mu)(K_\mu(w))$. Then if $\nu$ is positive, the
deformation in the direction of $\psi$ is tangent to a curve of positive measures. Note
also that we have the inverse formula,
\begin{equation*}
G_\nu(w) = G_\mu(w) + \psi(G_\mu(w)) \cdot G_\mu'(w)
\end{equation*}
\end{enumerate}
\end{Remark}

\subsection{Discussion}
The theorem can be interpreted as follows. Let $G$ be a univalent function on $\C^+$.
Define $K=G^{-1}$, $\Rtr(G)(z) = K(z) - 1/z$, and let $\psi$ be analytic on $G(\C^+)$.
Then the derivative of the map $\Rtr^{-1}$ (which is nothing other than the differential
of the operation of functional inversion) at a point $\Rtr(G)$ in the direction $\psi$ is
$- G' \psi(G)$. This linear map is invertible; the inverse linear map (which is the
differential of $\Rtr$) is $\psi \mapsto -K' \psi(K)$ (where $\psi$ is now analytic in
$\C^+$). Finally, let $\mathcal{T}_\alpha^\Rtr (\psi)(z) = 2\beta \psi(\beta z)$ and
denote $G_\alpha := G_{\varphi_\alpha}$. Then the derivative of $\mathcal{T}_\alpha =
\Rtr^{-1} \circ \mathcal{T}_\alpha^\Rtr \circ \Rtr$ at $G_\alpha$ is
\begin{equation*}
\begin{split}
D_{G_\alpha} \mathcal{T}_\alpha (\psi)(z)&= 2G_\alpha'(z) \beta K_\alpha'(\beta
G_\alpha(z)) \psi(K_\alpha(\beta G_\alpha(z))) \\
 &= 2 \psi(\omega_\alpha(z)) \omega_\alpha'(z)
\end{split}
\end{equation*}
Here
\begin{equation*}
\begin{split}
\omega_\alpha(z) &= K_{\varphi_\alpha}(\beta G_{\varphi_\alpha}(z)) =
\frac{1}{\beta} \omega_{\varphi_\alpha \circ S_{\beta}, \varphi_\alpha } \\
&= \frac{1}{\beta} \omega_{\varphi_\alpha \circ S_{\beta},
\varphi_\alpha \circ S_{\beta} \boxplus \varphi_\alpha \circ S_{\beta}} =
\frac{1}{\beta} K_{\varphi_\alpha \circ S_{\beta}} \circ G_{\varphi_\alpha \circ S_{\beta}
\boxplus \varphi_\alpha \circ S_{\beta}}
\end{split}
\end{equation*}
is a particular instance of the transition probability function of \cite{Voi93,Bia95}.

The eigenfunctions of $\mathcal{T}_\alpha^\Rtr$ on the space of all analytic functions in
the upper half plane are of the form $t e^{i\phi} z^a$, with eigenvalues
$2^{1-(a+1)/\alpha}$. Restricting to various spaces selects particular values of $\phi,
a$. Thus the eigenfunctions for the differential of $\mathcal{T}_\alpha$ (resp.,
$T_\alpha$) are the (boundary values of) the functions $e^{i \phi} G_\alpha' G_\alpha^a$.

\begin{Ex}
For the free normal (semicircular) case $\alpha=2$, $\nu = \chi$ the Cauchy transforms of
the eigenfunctions of the operator $D T$ (which are the eigenfunctions of the operator $D
\mathcal{T}$) are given by
\begin{equation*}
G_\chi' G_\chi^a = e^{i \phi} \frac{1}{\sqrt{z^2 - 4}} \left( \frac{z - \sqrt{z^2 - 4}}{2}
\right)^{x+yi}
\end{equation*}
By taking the boundary values $- \frac{1}{\pi} \im G(t + 0i)$ (see \cite{Akh65,Hor90}),
the eigenfunctions themselves are
\begin{multline*}
\frac{1}{\sqrt{t^2 -4}} {\abs{ \frac{t - \sqrt{t^2 - 4}}{2} }}^x e^{- y \pi}
\sin \left( y \log \abs{ \frac{t - \sqrt{t^2 -4}}{2} } + \phi + x \pi \right)
\chf_{(- \infty, -2]} (t) \,dt\\
+\frac{1}{\sqrt{t^2 -4}} {\abs{ \frac{t - \sqrt{t^2 - 4}}{2} }}^x e^{- y \pi}
\sin \left( y \log \abs{ \frac{t - \sqrt{t^2 -4}}{2} } + \phi  \right)
\chf_{[2, \infty)} (t) \,dt
 \\
+\frac{1}{\sqrt{4 - t^2}} \exp \left( y \cos^{-1}(t/2) \right)
\cos \left( x \cos^{-1}(t/2) - \phi \right) \chf_{[-2, 2]}(t) \,dt
\end{multline*}
It is easy to see that the Criterion~\ref{nec} requires that (for some $\phi$) we have $x
= \re a \geq 1$ or ($x= \re a \geq -1, y = \im a =0$); the corresponding point spectrum is
the union of the unit disc and the interval $[1,2]$. On the other hand, for $-1 < x < 1,
y=0$ the functions are the $R$-transforms of freely stable distributions, while for $x
\in \N,\, y=0$ the functions are entire. Thus by Criterion~\ref{suf}\ref{deform} the
corresponding eigenfunctions are in the tangent space to the space of positive measures.

All the moments of a measure are finite iff its Cauchy transform is analytic at infinity.
If a function $G$ is defined by the above expression on $\C^+$ and satisfies $G(\bar{z}) =
\overline{G(z)}$, it is analytic at infinity iff $a \in \N$ and $\phi=0$. Notice that
among the above measures, these are precisely the compactly supported ones. Explicitly
their Cauchy transforms are $- \frac{1}{\sqrt{z^2 - 4}} \left( \frac{z - \sqrt{z^2 -
4}}{2} \right)^n$. The eigenfunctions themselves are
\begin{equation*}
\frac{1}{\sqrt{4 - x^2}} \cos(n \cos^{-1}(x/2)) \chf_{[-2,2]}(x) \,dx
\end{equation*}
That is, we recover the Chebyshev functions of the first kind.
\end{Ex}

\begin{Remark}
Here we can see another difference from the classical case. As noted above, by a result of
Bercovici and Voiculescu \cite{BV95} the deformations in the directions $z^n, n \in \N$
actually produce positive measures (for small enough time). This is in contrast with a
classical theorem of Marcinkiewicz, which states that for $P$ a polynomial, $e^P$ is never
a characteristic function (i.e. a Fourier transform of a positive measure) if the degree of
$P$ is greater than 2 (see e.g. \cite[Thm 3.13]{Ram67}).
\end{Remark}

\begin{Ex}
For the 1-stable symmetric distribution, which is the Cauchy distribution $\varphi_1$, the
eigenfunction Cauchy transforms are $\frac{1}{(z-i)^a}$, and the eigenfunctions are $D^a
\frac{1}{x^2 + 1} \,dx$. Notice that these are exactly the same as in the classical case. Therefore not only are
the free 1-stable distributions the same as classical ones \cite{BV93}, but their small
neighborhoods look the same as well.
\end{Ex}

\begin{Remark}
For a general freely stable distribution, using for example the formula $D_G
\Rtr^{-1}(\psi) = - \frac{\psi}{K'}(G)$, one can obtain parametric expression \`{a} la
Biane \cite[Appendix]{BPB96} for the densities of the corresponding boundary values. In
particular, one has such expressions for the densities of the eigenfunctions of the free
stable central limit operators. It is not clear whether they are of use.
\end{Remark}

\subsection{Composition operators}
Finally, we have a brief discussion of the connections with the theory of composition
operators (see e.g. \cite{Val31,CM95}). The action of the operator $\mathcal{T}_\alpha$ on
the primitive (in $\C^+$) of a Cauchy transform is, up to a multiplicative constant 2 and
up to an additive constant, just the composition with the function $\omega_\alpha$.

One of the main theorems about composition operators is that such an operator is
necessarily conjugate to an operator of composition with a linear function \cite[Thm.~
2.53]{CM95}. In our case, in the terminology of \cite[Section 2.4]{CM95} the operator
$\mathcal{T}_\alpha$ has a natural halfplane-dilation model provided by conjugating with
the linearization of the $R$-transform, i.e. $G_\alpha \circ \omega_\alpha = \beta
G_\alpha$. A \emph{fundamental set} \cite[Defn.~2.54]{CM95} for $\psi$ is  an open,
connected, simply connected domain $\Delta$ such that $\psi(\Delta) \subset \Delta$ and the
iterates of any compact set end up in it after a finite number of steps. It is not hard to
see (e.g. \cite[Appendix]{BPB96}) that $\C^+$ serves as a fundamental set for both
$G_\alpha$ and $\omega_\alpha$ while $G(\C^+)$ is a fundamental set for $K_\alpha$.

A number of results on the spectra of composition operators on various classical spaces are
known. In particular, on a Hardy space $H^\infty$ the spectral radius of a composition
operator is equal to 1 \cite[3.1]{CM95}. In our case the composition operator is defined on
tangent spaces to a certain cone in $H^\infty$.

\end{document}